\documentclass[reqno,11pt]{article}
\usepackage{mathrsfs}
\usepackage{amssymb,amsmath}
\usepackage{amsthm}
\usepackage{amsmath}
\usepackage{amsfonts}
\usepackage{enumerate}
\usepackage{color}
\usepackage{bm}
\DeclareMathOperator{\X}{X}
\usepackage{graphicx}

\numberwithin{equation}{section}

\newtheorem{theorem}{Theorem}[section]
\newtheorem{lemma}[theorem]{Lemma}
\newtheorem{corollary}[theorem]{Corollary}
\newtheorem{proposition}[theorem]{Proposition}

\theoremstyle{definition}
\newtheorem{definition}[theorem]{Definition}
\newtheorem{assumption}[theorem]{Assumption}
\newtheorem{example}[theorem]{Example}

\theoremstyle{remark}
\newtheorem{remark}[theorem]{Remark}

\begin{document}
\title{Smooth Symmetric Transonic Isothermal Flows with Nonzero Angular Velocity}
\author {Na Zhang\thanks{School of mathematics and statistics, Wuhan University, Wuhan, Hubei Province, 430072, People's Republic of China. Email: nzhang@whu.edu.cn}}
\date{}
\maketitle

\def\be{\begin{eqnarray}}
\def\ee{\end{eqnarray}}
\def\be{\begin{equation}}
\def\ee{\end{equation}}
\def\ba{\begin{aligned}}
\def\ea{\end{aligned}}
\def\bay{\begin{array}}
\def\eay{\end{array}}
\def\bca{\begin{cases}}
\def\eca{\end{cases}}
\def\p{\partial}
\def\hphi{\hat{\phi}}
\def\bphi{\bar{\phi}}
\def\no{\nonumber}
\def\eps{\epsilon}
\def\de{\delta}
\def\De{\Delta}
\def\om{\omega}
\def\Om{\Omega}
\def\f{\frac}
\def\th{\theta}
\def\vth{\vartheta}
\def\la{\lambda}
\def\lab{\label}
\def\b{\bigg}
\def\var{\varphi}
\def\na{\nabla}
\def\ka{\kappa}
\def\al{\alpha}
\def\La{\Lambda}
\def\ga{\gamma}
\def\Ga{\Gamma}
\def\ti{\tilde}
\def\wti{\widetilde}
\def\wh{\widehat}
\def\ol{\overline}
\def\ul{\underline}
\def\Th{\Theta}
\def\si{\sigma}
\def\Si{\Sigma}
\def\oo{\infty}
\def\q{\quad}
\def\z{\zeta}
\def\co{\coloneqq}
\def\eqq{\eqqcolon}
\def\di{\displaystyle}
\def\bt{\begin{theorem}}
\def\et{\end{theorem}}
\def\bc{\begin{corollary}}
\def\ec{\end{corollary}}
\def\bl{\begin{lemma}}
\def\el{\end{lemma}}
\def\bp{\begin{proposition}}
\def\ep{\end{proposition}}
\def\br{\begin{remark}}
\def\er{\end{remark}}
\def\bd{\begin{definition}}
\def\ed{\end{definition}}
\def\bpf{\begin{proof}}
\def\epf{\end{proof}}
\def\bex{\begin{example}}
\def\eex{\end{example}}
\def\bq{\begin{question}}
\def\eq{\end{question}}
\def\bas{\begin{assumption}}
\def\eas{\end{assumption}}
\def\ber{\begin{exercise}}
\def\eer{\end{exercise}}
\def\mb{\mathbb}
\def\mbR{\mb{R}}
\def\mbZ{\mb{Z}}
\def\mc{\mathcal}
\def\mcS{\mc{S}}
\def\ms{\mathscr}
\def\lan{\langle}
\def\ran{\rangle}
\def\lb{\llbracket}
\def\rb{\rrbracket}
\def\fr#1#2{{\frac{#1}{#2}}}
\def\dfr#1#2{{\dfrac{#1}{#2}}}
\def\u{{\textbf u}}
\def\v{{\textbf v}}
\def\w{{\textbf w}}
\def\d{{\textbf d}}
\def\nn{{\textbf n}}
\def\x{{\textbf x}}
\def\e{{\textbf e}}
\def\U{{\textbf U}}
\def\M{{\textbf M}}
\def\FF{{\mathcal F}}
\def\X{{\mathcal X}}
\def\ww{{\mathcal w}}
\def\I{{\mathcal I}}
\def\S{{\mathcal S}}
\def\P{{\mathcal P}}
\def\W{{\mathcal W}}
\def\div{{\rm div\,}}
\def\curl{{\rm curl\,}}
\def\R{{\mathbb R}}
\def\N{{\mathbb N}}
\def\D{{\mathbb D}}
\def\O{{\mathbb O}}
\def\T{{\mathbb T}}
\def\g{{\textbf g}}
\def\F{{\textbf F}}
\def\A{{\textbf A}}
\def\f{{\textbf f}}
\def\r{{\textbf r}}
\def\mE{\mathcal{E}}
\def\Div{{\rm div}}
\def\i{\mathfrak{i}}

\begin{abstract}
  In this paper, the steady inviscid flows with radial symmetry for the isothermal Euler system are studied in an annulus. We present a complete classification of transonic radially symmetric flow patterns in term of physical boundary conditions at the inner and outer circle. By solving the one side boundary problem, we obtain that there exist accelerating or decelerating smooth transonic flows in an annulus. Moreover, the structural stability of these smooth symmetric transonic flows with nonzero angular velocity are further investigated. Furthermore, we examine the transonic solutions with shocks as well via prescribing suitable boundary conditions on the inner and outer circle. \\
\end{abstract}

\begin{center}
\begin{minipage}{5.5in}


\end{minipage}
\end{center}

\section{Introduction}\noindent

This paper is to study the isothermal flow in an annulus $\Omega=\{(x_1, x_2):r_0<r=\sqrt{x_1^2+x_2^2}<r_1\}$, which is governed by the following steady compressible Euler system
\begin{equation}\label{EPS}
\begin{cases}
\p_{x_1}(\rho u_1)+\p_{x_2}(\rho u_2)=0,\\
\p_{x_1}(\rho u_1^2)+\p_{x_2}(\rho u_1u_2)+\p_{x_1}p=0,\\
\p_{x_1}(\rho u_1u_2)+\p_{x_2}(\rho u_2^2)+\p_{x_2}p=0,\\
\end{cases}
\end{equation}
where $\rho$, $\u=(u_1, u_2)^t$, $p$ stand for the density, velocity, pressure respectively. 
And particularly, the isothermal gas for which the pressure $p$ is given by
$$
p=A\rho,
$$
here $A$ is a positive constant. In addition, the Bernoulli's function is given by
\be\label{B}
B=\frac12|{\bf u}|^2+A\ln\rho.
\ee

It should be emphasized that the $sound\,\, speed$ in our work is described by $c=\sqrt{p_{\rho}'(\rho)}=\sqrt{A}$. Depending on the $Mach \,\,number$ $\M=\frac{|\bf u|}{\sqrt{A}}$, the analytic and physical feature of \eqref{EPS} vary. When $\M>1$, the corresponding state of flow is called $supersonic$; When $\M<1$, the corresponding state is called $subsonic$.

In \cite{Courant1948}, Courant and Friedrichs employed the Hodograph transformation to rewrite the Euler system into a linear second order partial differential equation on the plane of the flow speed and angle, and obtained some special flows without concerning the boundary conditions. Among these special flows, there are circulatory flows and purely radial flows, which are symmetric flows with only angular and radial velocity respectively. Then their superpositions are called spiral flows. They had shown that the spiral flows may appear only outside a limiting circle, at which the Jacobian of the hodograph transformation is zero, and the flows may change smoothly from subsonic to supersonic or vice verse. However, there are few results about transonic spiral flows on the physical plane. Motivated by this, the authors in \cite{Weng1} investigated the steady inviscid compressible Euler flows in an annulus, they elucidated the effects of the angular velocity on radially symmetric transonic spiral flows. They have discovered some new transonic flow patterns and analyzed their properties in detail. And the structural stability of these transonic flows were investigated later in \cite{Weng2}. This inspired our interests to expand these results to the two-dimensional steady isothermal Euler flows.

In this work, we target on studying radially symmetric transonic spiral flows in an annulus associated with suitable boundary conditions on the inner and outer circle. And then the structural stability will further to be surveyed. Furthermore, in the presence of shock, we will give a complete classification of wave patterns according to the exit pressure and the position of the outer circle. It is convenient to use the polar coordinates
\begin{align}
&(x_1, x_2)=(r\cos\theta, r\sin\theta),\\
&(u_1, u_2)=(U_1\cos\theta-U_2\sin\theta, U_1\sin\theta+U_2\cos\theta)
\end{align}
to rewrite \eqref{EPS} as
\begin{eqnarray}\label{ProblemImm}
\begin{cases}
\partial_r(\rho U_1)+\frac{1}{r}\partial_{\theta}(\rho U_2)+\frac{1}{r} \rho U_1=0,\\
(U_1\partial_r +\frac{U_2}{r}\partial_{\theta}) U_1+\frac{1}{\rho} \partial_r p-\frac{U_2^2}r=0,\\
(U_1\partial_r +\frac{U_2}{r}\partial_{\theta}) U_2+\frac{1}{r\rho} \partial_{\theta}p+\frac{U_1U_2}r=0.\\
\end{cases}
\end{eqnarray}

Following the definition by Bers \cite{Lipman2}, we call a sonic point in a $C^2$ transonic flow is exceptional if the velocity is orthogonal to the sonic curve at this point. Therefore all sonic points of smooth transonic spiral flows here are nonexceptional and noncharacteristically degenerate in the presence of the nonzero angular velocity. Wang and Xin proved the existence and uniqueness of the smooth transonic flow of Meyer type in de Laval nozzles in \cite{WX1, WX4, WX2, WX3}, and all sonic points are exceptional and characteristically degenerate.

In \cite{GC1, Xie1, Xie2}, the authors proved the existence of subsonic-sonic weak solutions to the 2-D steady potential equation through the compensated compactness. Later on, the subsonic-sonic limit for multidimensional potential flows and steady Euler flows are examined in \cite{GC2, FH}. Nevertheless, the solutions obtained by the sonic-sonic limit merely satisfy the equations in the sense of distribution, and there is no information about the regularity and degeneracy properties near sonic points and their distribution in flow region. In \cite{WX1, WX4, WX2, WX3}, the authors studied the subsonic-sonic flows by using the Chaplygin equations in the plane of the velocity potential and the stream function as well as the comparison principles. However, in the presence of the nonzero angular velocity, it is expected that the sonic points of the background flows are nonexceptional, and the transformed sonic curve in the potential-stream functions plane is not a straight line in general, which is different from the cases in \cite{WX1, WX4, WX2, WX3}. The approach adapted by Wang and Xin can not work in this case.

Recently, Weng, Xin and Yuan investigated the steady inviscid compressible Euler flows with radial symmetric in an annulus in the presence of the nonzero angular velocity in \cite{Weng1}. They gave a complete classification of flow patterns in terms of boundary conditions at the inner and outer circle. By considering the one side boundary value problem, they showed that there exist acclerating or decelerating smooth transonic flows  with all sonic points being nonexceptional and noncharacteristically degenerate. Furthermore, they examined the transonic flows with shocks. Due to the fact that the radial velocity will jump from supersonic to subsonic and the angular velocity experiences no jump across the shock, which leads to the total velocity after the shock may be supersonic, sonic and subsonic. The authors discovered that besides the well-known supersonic-subsonic shock solutions in a divergent nozzle as in the case without the angular velocity, there are supersonic-supersonic shock solutions and the flow at downstream may change smoothly from supersonic to subsonic. A new and interesting phenomenon picked up that there also exist a supersonic-sonic shock solution where the shock front and the sonic circle coincide. Later on, they further investigated the structural stability of this special transonic flows in \cite{Weng2}. They considered the structural stability of smooth cylindrically symmetric transonic flows in a concentric cylinder under the cylindrical and axi-symmetric perturbations. As we all know, the steady Euler system is elliptic-hyperbolic mixed in subsonic regions and degenerates at the sonic points. To overcome this difficulty, Weng, Xin and Yuan applied the deformation-curl decomposition for the steady Euler system established in \cite{Weng3, Weng4} to analyze smooth symmetric transonic flows with nonzero angular velocity and vorticity.

As far as the isothermal transonic spiral flows with radial symmetry are considered, we will continue to explore their properties in an annulus. We will expand and integrate these results established in \cite{Weng1} and \cite{Weng2} into the isothermal flows in this paper. In our work, we will also give a complete classification of flow patterns according to the exit pressure and the position of the outer or inner circle and prove there exist accelerating or decelerating smooth transonic isothermal flows, with the aim to further investigate their structural stability. Furthermore, the transonic isothermal flows with shocks are studied as well.

The rest of the paper will be organized as follows. In Section \ref{smooth}, we will investigate radially symmetric transonic spiral isothermal flows in an annulus by considering the one side boundary value problem to \eqref{EPS} and show that there are accelerating or decelerating smooth transonic flows. Section \ref{stability} is devoted to investigate the structural stability of smooth transonic spiral solutions, of which the key element is to seek an appropriate multiplier to this linearized mixed equations, which is of great importance to establish the basic and higher order energy estimates. In Section \ref{transonic}, in the presence of shocks, the precise description of transonic flows and new wave patterns will be presented by studying \eqref{EPS} with two sides boundary data.

\section{Smooth radially symmetric flows with nonzero angular velocity}\label{smooth}\noindent

It deserves to be mentioned that the steady two-dimensional isothermal flow in this paper is investigated in an annulus $\Omega=\{(x_1, x_2):r_0<r=\sqrt{x_1^2+x_2^2}<r_1\}$. In this section, we will study the radially symmetric smooth solution to \eqref{EPS} in such an annulus, where the fluid flows from the inner circle to the outer one and the outer circle to the inner one respectively.

For the radially symmetric solutions to \eqref{ProblemImm} in $\Omega$, the steady Euler system further reduces to
\begin{align}\label{radial-euler}
\begin{cases}
\frac{d}{dr}(\rho U_1)+\frac1r \rho U_1=0, \ \ &r_0<r<r_1,\\
\rho U_1U_1'-\frac{\rho U_2^2}{r}+A\rho'=0,\ \ &r_0<r<r_1,\\
U_1U_2'+\frac{U_1U_2}r=0,\ \ &r_0<r<r_1.
\end{cases}
\end{align}

\subsection{The flows moving from the inner to outer circle}\noindent

We now consider the first case: the fluid moves from the inner circle to the outer circle. The following boundary conditions will be posed on the entrance $r=r_0$:
\begin{eqnarray}\label{bd1}
\rho(r_0)=\rho_{0}>0, \ U_{1}(r_0)=U_{10}>0, \ U_{2}(r_0)=U_{20}.
\end{eqnarray}
For later use, we note that the boundary value problem \eqref{radial-euler} and \eqref{bd1} is equivalent to
\begin{align}\label{alg1}
\begin{cases}
\rho U_1=\dfr{m_1}{r},\ m_1= r_0 \rho_{0} U_{10},\\
U_2=\dfr{m_2}{r},\ m_2= r_0 U_{20},\\
B\equiv B_0, B_0=\dfrac{1}{2}(U_{10}^2+U_{20}^2)+A \ln\rho_{0}.\\
\end{cases}
\end{align}

Then the main result in this case is stated as follow.
\bt\label{ProbI}
The results are divided into three cases.
\begin{enumerate}[(i)]
  \item If the radial Mach number at $r=r_0$ satisfies $M_1^2(r_0)>1$, the boundary value problem \eqref{radial-euler} and \eqref{bd1} has a unique supersonic smooth solution in $\Omega$.
  \item If the total Mach number $|{\bf M}|^2(r_0)<1$, the problem \eqref{radial-euler}-\eqref{bd1} has a unique subsonic smooth solution in $\Omega$.
  \item If $|{\bf M}|^2(r_0)>1>M_1^2(r_0)$, there exists a constant $r_s>r_0$ defined in \eqref{critical-r}, which depends only on $r_0$ and the incoming flow, such that
  \begin{enumerate}[(1)]
  \item if $r_1>r_s$, there exists a unique smooth transonic solution to \eqref{radial-euler}-\eqref{bd1} in $\Omega$;
  \item if $r_1<r_s$, there exists a unique smooth supersonic solution to \eqref{radial-euler}-\eqref{bd1} in $\Omega$;
  \item if $r_1=r_s$, there exists a unique smooth supersonic-sonic solution to \eqref{radial-euler}-\eqref{bd1} in $\Omega$ with sonic circle located at $r=r_s$.
  \end{enumerate}
\end{enumerate}
\et

\bpf
For convenience, we set
\be
K(r)=B-\frac{1}{2}\fr{m_2^2}{r^2}.
\ee
Direct calculation gives the following expressions of derivative
\begin{align}
\begin{cases}\label{daoshu1}
\frac{\rho'}{\rho}=\frac{U_1^2+U_2^2}{r(A-U_1^2)}=\dfr{|\M|^2}{r(1-M_{1}^2)},\\
\frac{U_1'}{U_1}=-\frac{U_1^2+U_2^2}{r(A-U_1^2)}-\frac1 r=-\frac{A+U_2^2}{r(A-U_1^2)}=-\dfr{1+ M_{2}^2}{r(1-M_{1}^2)},\\
U_2'=-\frac{U_2}{r},\\
(M_1^2)'=-\frac{U_1^2}{r(A-U_1^2)}(2+2M_2^2)=\frac{M_1^2}{r(M_1^2-1)}(2+2M_2^2),\\
(M_2^2)'=-\frac{2U_2^2}{r A}=-\frac{2M_2^2}{r},\\
(|\M|^2)'=-\frac{2|\M|^2}{r(1-M_1^2)},
\end{cases}
\end{align}
where $M_i^2:=\frac{U_i^2}{A}$ for $i=1,2$, and ${\bf M}=(M_1,M_2)^t$.
\begin{enumerate}[(1)]
  \item Suppose that $M_1^2(r_0)-1>0$. It follows from $\eqref{daoshu1}_4$ that $M_1^2(r)-1$ increases as $r$ increases. Thus $|{\bf M}|^2(r)>M_1^2(r)>1$ for any $r>r_0$ and the flow is always supersonic. The existence and uniqueness of the solution follow directly from the standard theory of ODE system.

  \item Suppose that $M_1^2(r_0)<|\M|^2(r_0)<1$. Then by $\eqref{daoshu1}_4$, $\frac{d}{dr} M_1^2(r)<0$ and thus $M_1^2(r)<1$ for any $r>r_0$. This together with $\eqref{daoshu1}_6$ and $|\M|^2(r_0)<1$, implies that $|\M|^2(r)<1$ for any $r>r_0$ and the flow is always subsonic.

  \item Suppose that $M_1^2(r_0)<1<|\M|^2(r_0)$. As shown in (2), one derives a unique smooth solution with $M_1^2(r)<1$ for any $r>r_0$.
  In order to analyze the properities of solutions, we intend to determine the position of the sonic circle as follows. Denote the sonic circle by $r_s$, simple computation yields the value of the corresponding critical density
  \begin{align}\label{rhoc}
  \rho_s=\exp\left(\frac{B_0}{A}-\frac1 2\right).
  \end{align}
  Since $U_2$ is obtained and $U_1=\frac{m_1}{r\rho}$, we see that $\rho$ is a root of the following algebraic equation
  \begin{align}\label{Fr}
  Q_r(\rho):=A\rho^2 \ln\rho-(B_0-\dfr12
  \dfr{m_2^2}{r^2})\rho^2+\dfr{m_1^2}{2r^2}=0.
  \end{align}
  By the direct calculation, the sonic circle $r=r_s$, a root of $Q_r(\rho_s)$, can be derived in the form that
   \begin{eqnarray}\label{critical-r}
  r_s=\sqrt{\frac{m_1^2+ m_2^2\rho_s^2}{2(B_0-A\ln\rho_s)\rho_s^2}}.
  \end{eqnarray}
  Therefore if $r_1>r_s$ we obtain a transonic smooth solution on $\Omega$; if $r_1=r_s$ there exists a supersonic-sonic flow; if $r_1<r_s$, we obtain a purely hyperbolic flow.
\end{enumerate}

\epf

\subsection{The flows moving from the outer to inner circle}\label{subsec2}\noindent

Next, we consider another case, in which the fluid moves from the outer circle to the inner circle. The boundary conditions at the outer circle is given as
\begin{eqnarray}\label{bd2}
\rho(r_1)=\rho_{0}>0,\ U_{1}(r_1)=U_{10}<0, \ U_{2}(r_1)=U_{20}.
\end{eqnarray}

The main results for the boundary value problem \eqref{radial-euler} and \eqref{bd2} can be stated as follows.

\bt\label{ProbII}
{\it For a given incoming flow at the outer circle $r=r_1$, there exists a limiting circle $r^\sharp<r_1$ depending only on $r_1$ and the incoming flow state at $r_1$ such that the flow develops singularities when $r\to r^{\sharp}$ in the sense that if $r_0<r^\sharp$ there does not exist any global $C^1$-smooth solution to Problem of \eqref{radial-euler} and \eqref{bd2} in $\Omega$. If $r_0\geq r^\sharp$, there exists a smooth solution to \eqref{radial-euler} and \eqref{bd2} in $\Omega$, the flow pattern in $\Omega$ can be divided into the following cases.
\begin{enumerate}[(i)]
  \item Suppose that $M_1^2(r_1)>1$ or $|{\bf M}|^2(r_1)>1>M_1^2(r_1)$, then there exists an unique smooth supersonic solution to \eqref{radial-euler} and \eqref{bd2} in $\Omega$.
  \item Suppose that $|{\bf M}|^2(r_1)<1$, there exists a constant $0<r=r_s<r_1$ defined in \eqref{critical-r}, which depends only $r_1$ and the incoming flow, such that
      \begin{enumerate}[(1)]
      \item if $r_0>r_s$, there exists a unique smooth subsonic solution to \eqref{radial-euler} and \eqref{bd2} in $\Omega$;
      \item if $r_0<r_s$, there exists a unique smooth transonic solution to \eqref{radial-euler} and \eqref{bd2} in $\Omega$;
      \item if $r_0=r_s$, there exists an unique smooth subsonic-sonic solution to \eqref{radial-euler} and \eqref{bd2} in $\Omega$ with sonic circle located at $r=r_s$.
      \end{enumerate}
\end{enumerate}
}\et

Similarly,  for smooth solutions, \eqref{alg1} is also available to the boundary value problem \eqref{radial-euler} and \eqref{bd2}, in which $m_1= r_1 \rho_{10} U_{10}$ and $m_2= r_1 U_{20}$.
\bl\label{Lemma2}
There exists a $0<r^\sharp< r_1$ such that
\begin{enumerate}[(1)]
  \item for any $r\in(r^\sharp,r_1)$, the algebraic system \eqref{alg1} has exactly two solutions $(\rho^-(r),U_1^-(r),U_2(r))$ and $(\rho^+(r),U_1^+(r),U_2(r))$ with $(M_1^-)^2(r)>1$ and $(M_1^+)^2(r)<1$.
  \item for $r=r^\sharp$, the algebraic system \eqref{alg1} has exactly one solution $(\rho^-(r^\sharp),U_1^-(r^\sharp),U_2(r^\sharp))=(\rho^+(r^\sharp),U_1^+(r^\sharp),U_2(r^\sharp))$ with $(M_1^-)^2(r^\sharp)=1=(M_1^+)^2(r^\sharp)$.
  \item for any $r\in(0,r^\sharp)$, the algebraic system \eqref{alg1} has no solution.
\end{enumerate}
\el
\bpf
Note that $\rho$ satisfies the algebraic equation $Q_r(\rho)=0$ with $Q_r(\rho)$ defined by \eqref{Fr}. Set $\tilde{r}=r_1 U_{20}/\sqrt{2B_0}$. For each fixed $0<r\leq \tilde{r}$, $Q_r(\rho)$ is monotonically increasing for $\rho\in(0, +\infty)$, and $\lim\limits_{\rho\rightarrow0}Q_r(\rho)=m_1^2/(2r^2)>0$. Thus $Q_r(\rho)> Q_r(0)>0$, $\forall$ $\rho\in(0, +\infty)$, then $Q_r(\rho)=0$ has no solution in $[0,+\infty)$ and the cases(3) is proved.

For each $r>\tilde{r}$, a direct computation gives
\begin{align}
\frac{d Q_r(\rho)}{d\rho}=2A\rho\ln\rho+A\rho-2(B_0-\fr12 \fr{m_2^2}{r^2})\rho.
\end{align}
Thus $Q_r(\rho)$ will first decrease as $\rho$ increases from $0$ to $\rho_*(r)$ and then increases when $\rho$ increases further, where
\begin{align}\label{rhostard}
\rho_*=\exp\bigg(\fr{K(r)}{A}-\frac1 2\bigg)<\rho_s.
\end{align}

Straightforward calculations yield
\begin{align*}
&Q_r(\rho_*(r))=A\ln\rho_*(r)\cdot\rho_*^2(r)-(B_0-\frac1 2\frac{m_2^2}{r^2})\rho_*^2(r)+\frac{m_1^2}{2r^2}\\
&\quad\quad\quad\quad=\left(K(r)-\frac A 2-B_0+\frac1 2\frac{m_2^2}{r^2}\right)\cdot\exp\left(\frac{2K(r)}{A}-1\right)+\frac{m_1^2}{r^2}\\
&\quad\quad\quad\quad=-\frac A 2\cdot\exp\left(\frac{2K(r)}{A}-1\right)+\frac{m_1^2}{r^2},\\
&Q_{\tilde{r}}(\rho_*(\tilde{r}))=\frac{m_1^2}{2\tilde{r}^2}>0,\\
&Q_{r_1}(\rho_*(r_{1}))=-\frac{A}{2}\exp\left(\frac{2K(r_1)}{A}-1\right)+\frac{m_1^2}{2r_1^2}\\
&\quad\quad\quad\quad\quad=-\frac A 2 \exp\left(2\ln\rho_0+\frac{U_{10}^2}{A}-1\right)+\frac1 2\rho_0^2U_{10}^2\\
&\quad\quad\quad\quad\quad=\frac A 2\left(\rho_0^2\cdot\frac{U_{10}^2}{A}-\exp(\ln\rho_0^2+\frac{U_{10}^2}{A}-1)\right)\\
&\quad\quad\quad\quad\quad=\frac A 2\rho_0^2\left(\frac{U_{10}^2}{A}-\exp(\frac{U_{10}^2}{A}-1)\right).
\end{align*}
Let $l(x)=e x-e^x$, then $l(0)=0$ and $l'(x)<0$ for any $x\in(0, 1)\cup(1, +\infty)$. This implies that $Q_{r_1}(\rho_*(r_{1}))=\frac A 2\rho_0^2 e^{-1}\cdot l\left(\frac{U_{10}^2}{A}\right)<0$ provided that $\frac{U_{10}^2}{A}\neq1$. Furthermore, $Q_r(\rho_*(r))$ is monotonically decreasing as a function of $r$ due to
\begin{align*}
\fr{d}{dr}Q_r(\rho_*(r))=-\exp\left(\frac{2K(r)}{A}-1\right)\frac{m_2^2}{r^3}-\frac{2m_1^2}{r^{3}}<0, \, \, \forall r>\tilde{r}.
\end{align*}
Hence there exists a unique $r^\sharp\in (\tilde{r},r_1)$ such that $Q_{r^\sharp}(\rho_*(r^{\sharp}))=0$. And if $r\in(\tilde{r},r^\sharp)$, $Q_r(\rho_*(r))>0$; if $r\in(r^\sharp,r_1)$, then $Q_r(\rho_*(r))<0$. Note that for each $r\in [r^{\sharp}, r_1)$, $Q_r(\rho)$ decreases first as $\rho$ increases from $0$ to $\rho_*(r)$ and then increases when $\rho$ increases further. Therefore for each $r\in(r^\sharp,r_1)$, $Q_r(\rho)=0$ has exactly two solution $\rho^\pm(r)$ with $0<\rho^-(r)<\rho_*(r)<\rho^+(r)$. For $r=r^\sharp$, $Q_r(\rho)=0$ has exactly one solution $\rho=\rho_*(r)$. For each $r\in (\tilde{r}, r^{\sharp})$, $Q_r(\rho)=0$ has no solution in $[0,+\infty)$. Furthermore, $\rho^-(r)<\rho_*(r)<\rho^+(r)$ is equivalent to $M_1^+(r)<1<M_1^-(r)$.

\epf

As stated in \cite{Weng1}, one can establish Theorem \ref{ProbII} according to Lemma \ref{Lemma2} and a similar argument as in the proof of Theorem \ref{ProbI}.

\section{The Structural stability of smooth transonic spiral flows}\label{stability}\noindent

In this section, we will survey the structural stability of the background flow to \eqref{ProblemImm} in a annulus $\Omega=\{(r, \theta): r_0<r=\sqrt{x_1^2+x_2^2}<r_1, \theta\in[0, 2\pi]\}$. In other words, we are searching for the solutions with the form $(U_1(r, \theta), U_2(r, \theta), \rho(r, \theta), p(r, \theta))$, which satisfies
\be\label{Euler3}
\begin{cases}
\p_r(\rho U_1)+\frac1 r\p_{\theta}(\rho U_2)+\frac1 r\rho U_1=0,\\
(U_1\p_r+\frac{U_2}{r}\p_{\theta})U_1+\frac{1}{\rho}\p_rp-\frac{U_2^2}{r}=0,\\
(U_1\p_r+\frac{U_2}{r}\p_{\theta})U_2+\frac{1}{r\rho}\p_{\theta}p+\frac{U_1U_2}{r}=0,
\end{cases}
\ee
with
\be\label{bd3}
\begin{cases}
U_1(r_0, \theta)-w_0U_2(r_0, \theta)=a_0+\eps g_0(\theta),\\
U_2(r_1, \theta)=a_1+\eps g_1(\theta),
\end{cases}
\ee
where $g_0$, $g_1$ are given periodic functions with period $2\pi$. $w_0$ is a constant which to be specified later and
\be
a_0=U_{b1}(r_0)-w_0U_{b2}(r_0), \quad\quad a_1=U_{20}.
\ee

To start with, we consider the potential flow. It follows from the vorticity being free that
\be
\frac1 r\p_r(rU_2)-\frac1 r\p_{\theta}U_1=0.
\ee
And the Bernoulli' laws implies
\be\label{32}
\rho=\exp\Bigg\{\frac1 A(B_0-\frac1 2(U_1^2+U_2^2))\Bigg\},
\ee
which gives that
\be\label{der} \begin{split}
&\p_r\rho=\frac{\rho}{A}(-U_1\p_rU_1-U_2\p_rU_2),\\
&\p_{\theta}\rho=\frac{\rho}{A}(-U_1\p_{\theta}U_1-U_2\p_{\theta}U_2),\\
&\p_rp=\rho(-U_1\p_rU_1-U_2\p_rU_2).
\end{split} \ee
Then \eqref{Euler3} can be reduced to the following boundary value problem:
\be\label{Euler5}
\begin{cases}
\p_r(r\rho U_1)+\p_{\theta}(\rho U_2)=0,\\
\frac1 r\p_r(rU_2)-\frac1 r\p_{\theta}U_1=0,\\
U_1(r_0, \theta)-w_0U_2(r_0, \theta)=a_0+\eps g_0(\theta),\\
U_2(r_1, \theta)=a_1+\eps g_1(\theta).
\end{cases}
\ee

Based on the background solution obtained in Theorem \ref{ProbII}, we have the following structural stability result.

\bt\label{irro}
Denote the Mach numbers by
$$
M_{b1}(r)=\frac{U_{b1}(r)}{A},\, M_{b2}(r)=\frac{U_{b2}(r)}{A},\, {\bf M}_b(r)=(M_{b1}(r), M_{b2}(r))^t,
$$
and let the background flow be given such that
\be\label{111}
1<|{\bf M}_b(r_0)|^2<2.
\ee
Assume that $g_0, g_1\in H^3(\T_{2\pi})$. Then for any constant $w_0$ with
\be\label{112}
w_0\notin\left(\frac{M_{b1}(r_0)M_{b2}(r_0)-\sqrt{|{\bf M}_{b}(r_0)|^2-1}}{1-M_{b1}^2(r_0)}, \frac{M_{b1}(r_0)M_{b2}(r_0)+\sqrt{|{\bf M}_{b}(r_0)|^2-1}}{1-M_{b1}^2(r_0)}\right),
\ee
there exists a small constant $\eps_0$ depending on the background flow, $w_0$, $g_0$ and $g_1$, such that for any $0<\eps<\eps_0$, the problem \eqref{Euler5} has a unique smooth transonic irrotational solution $(U_1, U_2)\in H^3(\Omega)\subset C^{1, \alpha}(\bar{\Omega})$ with the following estimate
\be\label{34}
||U_1-U_{b1}||_3+||U_2-U_{b2}||_3\leq C\eps.
\ee

Moreover, all the sonic points form a closed arc with a parametric representation $r=\eta(\theta)\in C^1(\T_{2\pi})$, \,$\forall$\, $\alpha\in(0, 1)$. The sonic surface is closed to the background sonic circle in the sense that
\be\label{35}
||\eta(\theta)-r_s||_{C^1(\T_{2\pi})}\leq C\eps.
\ee
\et

\begin{remark}\nonumber
{\it All the sonic points to the transonic irrotational solution obtained in Theorem \ref{irro} are nonexceptional and noncharacteristically degenerate.
}
\end{remark}

\begin{remark}\nonumber
{\it The radially symmetric smooth flows where the fluid moves from the inner to the outer circle are also structurally stable under the same perturbations as in \eqref{Euler5} within the class of irrotational flows.
}
\end{remark}

To get rid of the trouble that the potential function corresponding to the background flow defined as $\phi_b(r, \theta)=\int_{r_1}^rU_{b1}(\eta)d\eta+r_1U_{20}\cdot\theta$ is not periodic in $\theta$, which is resulted from the fact that $\Omega$ is not simple connected and the background flow has a nonzero circulation, we define the deviation between the flow and the background flow by
\be\label{36}
\hat{U}_1=U_1-U_{b1}, \quad \hat{U}_2=U_2-U_{b2}, \quad \hat{\rho}=\rho-\rho_{b},
\ee
it is clearly that $\hat{\U}, \hat{\rho}$ satisfy
\be\label{37}
\begin{cases}
\p_r(r(\rho_b\hat{U_1}+(U_{b1}+\hat{U}_1)\hat{\rho})+\p_{\theta}(\rho_b\hat{U_2}+(U_{b2}+\hat{U}_2)\hat{\rho})=0,\\
\frac1 r\p_r(r\hat{U}_2)-\frac1 r\p_{\theta}\hat{U}_1=0,\\
\hat{U}_1(r_0, \theta)-w_0\hat{U}_2(r_0, \theta)=\eps g_0(\theta),\\
\hat{U}_2(r_1, \theta)=\eps g_1(\theta).
\end{cases}
\ee

Define the potential function
$$
\phi(r, \theta)=\int_{r_1}^r\hat{U}_1(\tau, \theta)\tau+\int_0^{\theta}(r_1\hat{U}_2(r_1, \tau)+d_0)d\tau=\int_{r_1}^r\hat{U}_1(\tau, \theta)\tau+\int_0^{\theta}(\eps r_1g_1(\tau)+d_0)d\tau,
$$
note that $d_0$ is introduced here in order to guarantee the periodic of $\phi$ in $\theta$. Indeed, $d_0=-\eps r_1\bar{g}_1$ with $\bar{g}_1=\frac{1}{2\pi}\int_0^{2\pi}g_1(\tau)d\tau$, hence $\phi$ is periodic in $\theta$ with period $2\pi$, and satisfies
\be\label{38}
\p_r\phi=\hat{U}_1, \p_{\theta}\phi=r\hat{U}_2+d_0.
\ee

\subsection{Linearized problem}

We now concentrate on building the linearized problem. At first, upon recalling the derivatives of variables presented in \eqref{daoshu1}, one has
\begin{align}
\begin{cases}\label{daoshu}
\rho_b'(r)=\dfr{|\M_b|^2}{r(1-M_{b1}^2)}\rho_b, \quad\quad\quad\quad (M_{b1}^2)'=-\frac{M_{b1}^2}{r(1-M_{b1}^2)}(2+2M_{b2}^2),\\
U_{b1}'(r)=-\dfr{1+ M_{b2}^2}{r(1-M_{b1}^2)}U_{b1},\quad\quad (M_{b2}^2)'=-\frac{2M_{b2}^2}{r},\\
U_{b2}'(r)=-\frac{U_{b2}}{r},\quad\quad\quad\,\quad\quad\quad\quad (|\M_b|^2)'=-\frac{2|\M_b|^2}{r(1-M_{b1}^2)}.
\end{cases}
\end{align}
It follows from the first equation of \eqref{Euler5} that
\be\label{aa}
U_1\p_r\rho+\rho\p_r U_1+\frac{U_2}{r}\p_{\theta}\rho+\frac{\rho}{r}\p_{\theta} U_2+\frac1 r\rho U_1=0.
\ee
Straightforward calculation by utilizing applying \eqref{der} yields
\be\label{zzz}
\p_r{U}_1+\frac1 r\p_{\theta}U_2+\frac{U_1}{r}=\frac{U_1^2}{A}\p_r{U}_1+\frac{U_1U_2}{A}\p_r{U}_2+\frac{U_1U_2}{rA}\p_{\theta}U_1+\frac{U_2^2}{rA}\p_{\theta}U_2.
\ee
Substituting \eqref{36} into \eqref{zzz} gives
\begin{align*}
(A-U_1^2)\p_r^2\phi+&(A-U_1^2)\p_r{U}_{b1}+\frac A r(\hat{U_1}+U_{b1})+\frac{A-U_2^2}{r^2}\p_{\theta}^2\phi+\frac A r\p_{\theta}U_{b2}\\
&=\frac{2U_1U_2}{r}\p_{r\theta}^2\phi-\frac{U_1U_2}{r}\hat{U_2}+U_1U_2\p_{r}U_{b2}+\frac{U_1U_2}{r}\p_{\theta}U_{b1}+\frac{U_2^2}{r}\p_{\theta}U_{b2}\\
&=\frac{2U_1U_2}{r}\p_{r\theta}^2\phi-\frac{U_1U_2}{r}\hat{U_2}-\frac{U_1U_2}{r}U_{b2}+\frac{2U_1U_2}{r}\p_{\theta}U_{b1}+\frac{U_2^2}{r}\p_{\theta}U_{b2},
\end{align*}
that is to say,
\be\label{tt}
(A-U_1^2)\p_r^2\phi+\frac{A-U_2^2}{r^2}\p_{\theta}^2\phi-\frac{2U_1U_2}{r}\p_{r\theta}^2\phi=-(A-U_1^2)\p_r{U}_{b1}-\frac{A}{r}(\hat{U}_1+U_{b1})-\frac{U_1U_2}{r}U_2.
\ee
Continue to deal with the right side of \eqref{tt} by employing \eqref{daoshu}, we find it equals to
\begin{align*}
&\quad-A\p_r{U}_{b1}+(U_{b1}+\hat{U}_1)^2\p_r{U}_{b1}-\frac A r\p_r\phi-\frac A rU_{b1}-\frac1 r(U_{b1}+\hat{U}_1)(U_{b2}+\hat{U}_2)^2\\
&=\frac{A+U_{b2}^2}{r}U_{b1}-\frac{2(1+M_{b2}^2)}{r(1-M_{b1}^2)}U_{b1}^2\hat{U}_1-\frac{1+M_{b2}^2}{r(1-M_{b1}^2)}U_{b1}\cdot\hat{U}_1^2-\frac A r\p_r\phi-\frac{A}{r}U_{b1}-\frac1 r U_{b1}U_{b2}^2\\
&\quad-2U_{b1}U_{b2}\cdot\frac{\p_{\theta}\phi-d_0}{r^2}-\frac1 rU_{b1}\hat{U}_2^2-\frac 1 r\hat{U}_1U_2^2\\
&=-\frac{2(1+M_{b2}^2)}{r(1-M_{b1}^2)}U_{b1}^2\p_r\phi-\frac{1+M_{b2}^2}{r(1-M_{b1}^2)}U_{b1}\cdot\hat{U}_1^2
-\frac{2U_{b1}U_{b2}}{r^2}(\p_{\theta}\phi-d_0)\\
&\quad-\frac{U_{b1}}{r}\hat{U}_2^2-\frac{U_{b2}^2+A}{r}\hat{U}_1-\frac{1}{r}U_2^2\hat{U}_1+\frac{1}{r}U_{b2}^2\hat{U}_1.
\end{align*}
Then we derive the following second order mixed type equation
\begin{eqnarray}\label{mixed}
\begin{cases}
L\phi\equiv a_{11}\p_r^2\phi+a_{22}\p_{\theta}^2\phi+(a_{12}+a_{21})\p_{r\theta}^2\phi+e_1(r)\p_r\phi+e_2(r)\p_{\theta}\phi=f(U_1, U_2),\\
\p_r\phi(r_0, \theta)-w_0\frac{1}{r_0}\p_{\theta}\phi(r_0, \theta)=\epsilon g_0(\theta)-w_0\frac{1}{r_0}d_0,\\
\p_{\theta}\phi(r_1, \theta)=d_0+r_1\epsilon g_1(\theta),\\
\phi(r_1, 0)=0,
\end{cases}
\end{eqnarray}
where
\be
\begin{cases}
a_{11}(U_1, U_2)=A-U_1^2, \quad a_{22}(U_1, U_2)=\frac{A-U_2^2}{r^2}, \quad a_{12}(U_1, U_2)=a_{21}(U_1, U_2)=-\frac{U_1U_2}{r},\\
e_1(r)=\frac{2(1+M_{b2}^2)}{r(1-M_{b1}^2)}U_{b1}^2+\frac{U_{b2}^2+A}{r}, \quad e_2(r)=\frac{2U_{b1}U_{b2}}{r^2},\\
f(U_1, U_2)=e_2(r)d_0+U_{b1}'\cdot\hat{U}_1^2-\frac{U_{b1}}{r}\hat{U}_2^2-\frac{U_2^2-U_{b2}^2}{r}\hat{U}_1.
\end{cases}
\ee

Next we linearize the equation. First of all, representing the function space as
\be
\S=\{\phi\in H^4(\Omega):||\phi||_4\leq\delta_0\},
\ee
where $\delta_0>0$ will be specified later. Then for any function $\bar{\phi}\in\S$, define
$$
\bar{U}_1=U_{b1}+\p_r\bar{\phi}, \quad \bar{U}_2=U_{b2}+\frac1 r\p_{\theta}\bar{\phi}-\frac{d_0}{r}.
$$
We propose to construct an operator $\mathcal{M}: \bar{\phi}\in\S\mapsto \phi\in\S$, where $\phi$ is a solution of the following linear mixed type second-order partial differential equation
\begin{eqnarray}\label{linear mixed}
\begin{cases}
\bar{L}\phi\equiv a_{11}(\bar{U_1}, \bar{U_2})\p_r^2\phi+a_{22}(\bar{U_1}, \bar{U_2})\p_{\theta}^2\phi\\
\quad\quad\quad\quad+2a_{12}(\bar{U_1}, \bar{U_2})\p_{r\theta}^2\phi+e_1(r)\p_r\phi+e_2(r)\p_{\theta}\phi=f(U_1, U_2),\\
\p_r\phi(r_0, \theta)-w_0\frac{1}{r_0}\p_{\theta}\phi(r_0, \theta)=\epsilon g_0(\theta)-w_0\frac{1}{r_0}d_0=O(\epsilon),\\
\p_{\theta}\phi(r_1, \theta)=d_0+r_1\epsilon g_1(\theta)=O(\eps),\\
\phi(r_1, 0)=0.
\end{cases}
\end{eqnarray}
Let
$$
\xi(r)=\int_{r_0}^{r}\frac{M_{b1}M_{b2}(\tau)}{1-M_{b1}^2(\tau)}\frac{d\tau}{\tau},
$$
and
$$
y_1=r, \, y_2=\xi(r)+\theta,
$$
one has
\begin{align*}
&\p_r\phi=\p_{y_1}\phi\cdot \p_r y_1+\p_{y_2}\phi\cdot \p_r y_2=\p_{y_1}\hat{\phi}+{\xi}'(r)\p_{y_2}\hat{\phi},\\
&\p_{\theta}\phi=\p_{y_2}\phi\cdot \p_{\theta}y_2=\p_{y_2}\hat{\phi},\\
&\p_{r\theta}\phi=\p_{\theta}(\p_{y_1}\phi+\p_{y_2}\phi\cdot \xi'(r))=\p_{y_1y_2}^2\phi+\xi'(r)\p_{y_2}^2\phi,\\
&\p_{r}^2\phi=\p_r(\p_{y_1}\hat{\phi}+\p_{y_2}\hat{\phi}\xi'(r))=\p_{y_1}^2\hat{\phi}+2\p_{y_1y_2}^2\hat{\phi}\cdot \xi'(r)+\p_{y_2}^2\hat{\phi}(\xi'(r))^2+\p_{y_2}\hat{\phi}\cdot \xi''(r),\\
&\p_{\theta}^2\phi=\p_{\theta}(\p_{y_2}\phi\cdot \p_{\theta}y_2)=\p_{y_2}^2\phi\p_{\theta}y_2=\p_{y_2}^2\hat{\phi},
\end{align*}
by set the function $\hat\phi(y_1, y_2)=\phi(y_1, y_2-\xi(y_1))$.
Inserting these derivatives expressions into \eqref{linear mixed}, then the linear mixed-type second-order equation can be written as
\begin{eqnarray}\label{mix}
\begin{cases}
\hat{L}\hat{\phi}\equiv \sum\limits_{i,j=1}^{2}k_{ij}\p_{y_iy_j}^2\hat{\phi}+\sum\limits_{i=1}^2k_i\p_{y_i}\hat{\phi}=\hat{f}(U_1, U_2),\\
r_0\p_{y_1}\hat{\phi}(r_0, y_2)+(r_0\xi'(r_0)-w_0)\p_{y_2}\hat{\phi}(r_0, y_2)=g_2(y_2),\\
\p_{y_2}\hat{\phi}(r_1, y_2)=g_3(y_2),\\
\hat{\phi}(r_1, \xi(r_1))=0,
\end{cases}
\end{eqnarray}
where $(y_1, y_2)\in(r_0, r_1)\times\T_{2\pi}$ and
$$
k_{11}\equiv1, \quad k_{12}(\bar{U}_1, \bar{U}_2)=k_{21}(\bar{U}_1, \bar{U}_2)=\frac{a_{12}(\bar{U}_1, \bar{U}_2)+a_{11}(\bar{U}_1, \bar{U}_2)\xi'(y_1)}{a_{11}(\bar{U}_1, \bar{U}_2)},
$$
$$
k_{22}(\bar{U}_1, \bar{U}_2)=(\xi'(y_1))^2+\frac{a_{22}(\bar{U}_1, \bar{U}_2)+(a_{12}+a_{21})(\bar{U}_1, \bar{U}_2)\xi'(y_1)}{a_{11}(\bar{U}_1, \bar{U}_2)},
$$
$$
k_1(\bar{U}_1, \bar{U}_2)=\frac{e_1(y_1)}{a_{11}(\bar{U}_1, \bar{U}_2)}, \quad k_2(\bar{U}_1, \bar{U}_2)=\xi''(y_1)+\frac{e_1(y_1)\xi'(y_1)+e_2(y_1)}{a_{b11}(\bar{U}_1, \bar{U}_2)},
$$
$$
g_2(y_2)=r_0\epsilon g_0(y_2)-w_0d_0, \quad g_3(y_2)=d_0+r_1\epsilon g_1(y_2-\xi(r_1)).
$$
Elementary computations gives
$$
a_{b12}(y_1)+a_{b11}(y_1)\xi'(y_1)=0,
$$
which combines with the identity in proposition \ref{prop1}, $k_{b2}(y_1)=0$, we obtain the following estimates
\begin{equation}\label{estimate}
\begin{cases}
&||k_{12}(\bar{U}_1, \bar{U}_2)||_3+||k_{22}(\bar{U}_1, \bar{U}_2)-k_{b22}(y_1)||_3\leq C_0\delta_0,\\
&||k_{1}(\bar{U}_1, \bar{U}_2)-k_{b1}(y_1)||_3+||k_{2}(\bar{U}_1, \bar{U}_2)||_3\leq C_0\delta_0.
\end{cases}
\end{equation}

For ease of notation, we will always write $\hat{\phi}$ for $\phi$ in the rest of the paper.

\subsection{Some properties of the background flow}

In order to establish the basic estimate for the linearized mixed type potential equation, we now present some essential properties as follows. Let
\begin{align*}
&k_{b1}=\frac{e_1}{A-U_{b1}^2},\\
&k_{b2}(r)=\xi''(r)+\frac{e_1(r)\xi'(r)+e_2(r)}{A-U_{b1}^2},\\
&k_{b22}(r)=\frac{1-|\M_b(r)|^2}{r^2(1-M_{b1}^2(r))^2}.
\end{align*}
\bp\label{prop1}
Let $(U_{b1}, U_{b2}, \rho_b)$ be the background transonic flow, then we have the following identities for any $r\in[r_0, r_1]$:
\begin{align}
&k_{b2}\equiv0,\label{21}\\
&2k_{b1}k_{b22}+k'_{b22}(r)=\frac{2|\M_b|^2}{r^3(1-M_{b1}^2)^3}(2-|\M_b|^2).\label{22}
\end{align}
\ep
\bpf
Direct and explicit computations give
\be
\xi'(r)=\frac{M_{b1}M_{b2}(r)}{r(1-M_{b1}^2(r))},
\ee
and
\be \begin{split}
\xi''(r)&=\frac{(M_{b1}'M_{b2}+M_{b1}M_{b2}')(1-M_{b1}^2)\cdot r-M_{b1}M_{b2}\cdot(-2rM_{b1}M_{b1}'+1-M_{b1}^2)}{r^2(1-M_{b1}^2(r))^2}\no\\
&=\frac{M_{b1}'M_{b2}+M_{b1}M_{b2}'}{r(1-M_{b1}^2)}+\frac{M_{b1}M_{b2}}{r^2(1-M_{b1}^2(r))^2}\cdot(2rM_{b1}M_{b1}'-(1-M_{b1}^2))\no\\
&=\frac{-M_{b1}M_{b2}(1+M_{b2}^2)-2M_{b1}M_{b2}(1-M_{b1}^2)}{r^2(1-M_{b1}^2(r))^2}-\frac{2M_{b1}^3M_{b2}(1+M_{b2}^2)}{r^2(1-M_{b1}^2(r))^3}.
\end{split} \ee
Then we get
\begin{align*}
e_1(r)\xi'(r)+e_2(r)&=\frac{2(1+M_{b2}^2)}{r(1-M_{b1}^2(r))}+\frac{2U_{b1}U_{b2}}{r^2}+\frac{U_{b2}^2+A}{r}\cdot\frac{M_{b1}M_{b2}}{r(1-M_{b1}^2(r))}\no\\
&=\frac{2M_{b1}M_{b2}(1+M_{b2}^2)}{r^2(1-M_{b1}^2(r))^2}+\frac{2U_{b1}U_{b2}}{r^2}+\frac{M_{b1}M_{b2}(U_{b2}^2+A)}{r^2(1-M_{b1}^2(r))},
\end{align*}
which leads to
\begin{align*}
&\quad k_{b2}(r)=\xi''(r)+\frac{e_1(r)\xi'(r)+e_2(r)}{A-U_{b1}^2}\\
&=\frac{-M_{b1}M_{b2}(1+M_{b2}^2)}{r^2(1-M_{b1}^2(r))^2}-\frac{2M_{b1}^3M_{b2}(1+M_{b2}^2)}{r^2(1-M_{b1}^3(r))}
+\frac{2M_{b1}M_{b2}(1+M_{b2}^2)M_{b1}^2}{r^2(1-M_{b1}^3(r))}+\frac{M_{b1}M_{b2}(1+M_{b2}^2)}{r^2(1-M_{b1}^2(r))^2}\\
&=0.
\end{align*}
In additions, one can deduce that
\begin{align*}
k_{b22}'(r)&=\frac{(-|\M_b|^2)'}{r^2(1-M_{b1}^2)^2}-\frac{(1-|\M_b|^2)\cdot(2r(1-M_{b1}^2)^2+r^2\cdot2(1-M_{b1}^2)(-M_{b1}^2)')}{r^4(1-M_{b1}^2)^4}\\
&=\frac{2|\M_{b}|^2}{r^3(1-M_{b1}^2(r))^3}-\frac{2(1-|\M_b|^2)}{r^3(1-M_{b1}^2)^2}-\frac{2M_{b1}^2(1-|\M_b|^2)}{r^3(1-M_{b1}^2)^4}(2+2M_{b2}^2)\\
&=\frac{2|\M_{b}|^2}{r^3(1-M_{b1}^2(r))^3}-\frac{2(1-|\M_b|^2)(1-M_{b1}^2)}{r^3(1-M_{b1}^2)^3}
-\frac{4M_{b1}^2(1-|\M_b|^2)}{r^3(1-M_{b1}^2)^4}(1+|\M_b|^2-M_{b1}^2)\\
&=\frac{4|\M_{b}|^2-2-2M_{b1}^2+2|\M_b|^2M_{b1}^2}{r^3(1-M_{b1}^2(r))^3}-\frac{4|\M_{b}|^2M_{b1}^2(1-|\M_b|^2)}{r^3(1-M_{b1}^2(r))^4},
\end{align*}
and
as a consequence, we have
\begin{align*}
2k_{b1}k_{b22}+k_{b22}'(r)&=\frac{2(1-|\M_b|^2)}{r^2(1-M_{b1}^2)^2}\cdot\frac{(2(1+M_{b2}^2)M_{b1}^2}{r(1-M_{b1}^2)^2}
+\frac{2(1-|\M_b|^2)}{r^2(1-M_{b1}^2)^2}\cdot\frac{1+M_{b2}^2}{r(1-M_{b1}^2)}\\
&\quad+\frac{4|\M_{b}|^2-2-2M_{b1}^2+2|M_b|^2M_{b1}^2}{r^3(1-M_{b1}^2(r))^3}-\frac{4|\M_{b}|^2M_{b1}^2(1-|\M_b|^2)}{r^3(1-M_{b1}^2(r))^4}\\
&=\frac{4(1-|\M_b|^2)(1+M_{b2}^2)M_{b1}^2-4|\M_b|^2M_{b1}^2(1-|\M_b|^2)}{r^3(1-M_{b1}^2)^4}\\
&\quad+\frac{(2(1-|\M_b|^2)(1+M_{b2}^2)+4|\M_b|^2-2+2M_{b1}^2(|\M_b|^2-1)}{r^3(1-M_{b1}^2)^3}\\
&=\frac{2|\M_b|^2}{r^3(1-M_{b1}^2)^3}(2-|\M_b|^2).
\end{align*}

Therefore, the identities in the proposition are valid.
\epf

\subsection{Energy estimates for the linearized problem}

Based on \eqref{estimate}, we establish the basic and higher energy estimates to \eqref{mix} in this subsection provided that $k_{11}, k_{12}, k_{21}, k_{22}\in C^{\infty}(\bar{\Omega})$ and $g_2, g_3\in C^{\infty}(\T_{2\pi})$.

\bl\label{low}
Suppose that \eqref{111} and \eqref{112} hold, then there exists two constants $\sigma_*>0$, $\delta_*>0$ depending only on the background flow, and $m_0$, such that if $0<\delta_0\leq\delta_*$ in \eqref{estimate}, the solution to \eqref{mix} satisfies the following basic energy estimate
\be\label{goal}
||\phi||_1\leq\frac{C_*}{\sigma_*}(||\hat{f}||_{L^2(\Omega)}+\sum\limits_{j=0,1}||g_j||_{L^2(T_{2\pi})}),
\ee
where $C_*$ depends only on the $H^3(\Omega)$ norm of the coefficients $k_{ij}, k_i,$ for $i, j=1,2$.
\el
\bpf
We intend to derive the basic energy estimate, of which the key element is to find a multiplier. Let $w_1(y_1)$ and $w_2(y_1)$ be smooth functions of $y_1$ in $[r_0, r_1]$, which may be determined later. We first multiply the equation in \eqref{mix} by $w_1(y_1)\p_{y_1}\phi+w_2(y_1)\p_{y_2}\phi$ and integrate over $\Omega$, after integration by parts,
\begin{align*}
&\iint_{\Omega}\hat{f}(w_1(y_1)\p_{y_1}\phi+w_2(y_1)\p_{y_2}\phi)dy_1dy_2\\
&\quad=\iint_{\Omega}\p_{y_1}\left(\frac1 2(\p_{y_1}\phi)^2\cdot w_1\right)-\frac1 2(\p_{y_1}\phi)^2\cdot w_1'-\p_{y_2}(w_1k_{12})(\p_{y_1}\phi)^2\\
&\quad\quad\quad+\p_{y_2}(w_1k_{12}(\p_{y_1}\phi)^2)-\p_{y_1}(w_2k_{21})(\p_{y_2}\phi)^2+\p_{y_1}(w_2k_{21}(\p_{y_2}\phi)^2)+w_1\p_{y_1}\phi\cdot k_{22}\p_{y_2}^2\phi\\
&\quad\quad\quad+w_1k_1(\p_{y_1}\phi)^2+w_1k_2\p_{y_1}\phi\p_{y_2}\phi+\p_{y_1}(w_2\p_{y_1}\phi\p_{y_2}\phi)-\p_{y_1}(w_2\p_{y_2}\phi)\cdot\p_{y_1}\phi\\
&\quad\quad\quad+\p_{y_2}\left(\frac1 2(\p_{y_2}\phi)^2\cdot w_2k_{22}\right)-\p_{y_2}(w_2k_{22})\cdot\frac1 2(\p_{y_2}\phi)^2+w_2k_1\p_{y_1}\phi\p_{y_2}\phi+w_2k_2(\p_{y_2}\phi)^2dy_1dy_2\\
&\quad=\int_0^{2\pi}\frac1 2(\p_{y_1}\phi)^2\cdot w_1+w_2k_{21}(\p_{y_2}\phi)^2-\frac1 2w_1k_{22}(\p_{y_2}\phi)^2+w_2\cdot\p_{y_1}\phi\p_{y_2}\phi\Bigg|_{y_1=r_0}^{r_1}dy_2\\
&\quad\quad+\iint_{\Omega}-\frac1 2w_1'(\p_{y_1}\phi)^2-w_1\p_{y_2}k_{12}(\p_{y_1}\phi)^2-\p_{y_1}(w_2k_{21})(\p_{y_2}\phi)^2\\
&\quad\quad\quad\quad+\frac1 2\p_{y_1}(w_1k_{22})(\p_{y_2}\phi)^2-\p_{y_2}(w_1k_{22})\p_{y_1}\phi\p_{y_2}\phi+w_1k_1(\p_{y_1}\phi)^2+w_1k_2\p_{y_1}\phi\p_{y_2}\phi\\
&\quad\quad\quad\quad-w_2'\cdot\p_{y_1}\phi\p_{y_2}\phi-\frac1 2w_2\p_{y_2}k_{22}\cdot(\p_{y_2}\phi)^2+w_2k_1\p_{y_1}\phi\p_{y_2}\phi+w_2k_2(\p_{y_2}\phi)^2dy_1dy_2
\end{align*}
with the help of
\begin{align*}
w_1\p_{y_1}\phi\cdot k_{22}\p_{y_2}^2\phi&=\p_{y_2}(w_1k_{22}(\p_{y_1}\phi)^2)-\p_{y_2}(w_1k_{22}\p_{y_1}\phi)\p_{y_2}\phi\\
&=-w_1k_{22}\p_{y_1y_2}^2\phi\cdot \p_{y_2}\phi-\p_{y_2}(w_1k_{22})\p_{y_1}\phi\p_{y_2}\phi\\
&=\p_{y_1}\left(-\frac1 2w_1k_{22}(\p_{y_2}\phi)^2\right)+\frac1 2\p_{y_1}(w_1k_{22})(\p_{y_2}\phi)^2-\p_{y_2}(w_1k_{22})\p_{y_1}\phi\p_{y_2}\phi.
\end{align*}
Finally one obtains
\be\label{F} \begin{split}
&\quad\iint_{\Omega}\hat{f}(w_1(y_1)\p_{y_1}\phi+w_2(y_1)\p_{y_2}\phi)dy_1dy_2\\
&=\frac1 2\int_0^{2\pi}\Bigg(\left(\sqrt{w_1}\p_{y_1}\phi+\frac{w_2}{\sqrt{w_1}}\p_{y_2}\phi\right)^2+\left(-w_1k_{22}+2w_2k_{12}-\frac{w_2^2}{w_1}
\right)(\p_{y_2}\phi)^2\Bigg)dy_2\Bigg|_{y_1=r_0}^{r_1}\\
&\quad+\iint_{\Omega}\left(w_1k_1-\frac1 2w_1'-w_1\p_{y_2}k_{12}\right)(\p_{y_1}\phi)^2dy_1dy_2\\
&\quad+\iint_{\Omega}(w_2k_1-w_2'+w_1k_2-w_1\p_{y_2}k_{22})\p_{y_1}\phi\p_{y_2}\phi dy_1dy_2\\
&\quad+\iint_{\Omega}\left(\frac1 2\p_{y_1}(w_1k_{22})-\frac1 2w_2\p_{y_2}k_{22}-\p_{y_1}(w_2k_{12})+w_2k_2\right)(\p_{y_2}\phi)^2dy_1dy_2.
\end{split} \ee

We want next to derive an energy estimate, as discussed in \cite{Weng2}, it is necessary to show that there exist smooth functions $w_1(y_1)$ and $w_2(y_1)$ such that if $\delta_0$ is small enough, the following inequalities hold:
\be
\begin{cases}
w_1k_1-\frac1 2w_1'-w_1\p_{y_2}k_{12}\geq\sigma_*, \quad \forall(x, y)\in\Omega,\\
\frac1 2\p_{y_1}(w_1k_{22})-\frac1 2w_2\p_{y_2}k_{22}-\p_{y_1}(w_2k_{12})+w_2k_2\geq\sigma_*, \quad \forall(x, y)\in\Omega,\\
||w_2k_1-w_2'+w_1k_2-w_1\p_{y_2}k_{22}||_{L^{\infty}(\Omega)}\leq C_*\delta_0,\\
\left(k_{22}w_1+\frac{w_2^2}{w_1}-2k_{12}w_2\right)(r_0, y_2)>0, \quad \forall y_2\in \T_{2\pi},
\end{cases}
\ee
where $\sigma_*>0$ is a constant depending only on the background flow. Choose $w_1(y_1)$ such that
$$
w_1(y_1)k_{b1}(y_1)-\frac1 2w_1'(y_1)=\sigma_1,
$$
where $\sigma_1\in(0, 1)$ is a small positive constant to be determined later. Then one has
\be
w_1(y_1)=e^{\int_{r_0}^{y_1}2k_{b1}(\tau_1)d\tau_1}\left(w_1(r_0)-2\sigma_1\int_{r_0}^{y_1}e^{-\int_{r_0}^{y_1}2k_{b1}(\tau_2)d\tau_2}\right),
\ee
where $w_1(r_0)=1+2\int_{r_0}^{r_1}e^{-\int_{r_0}^{\tau_1}2k_{b1}(\tau_2)d\tau_2}d\tau_1.$ Then for any $y_1\in[r_0, r_1]$, one has $w_1(y_1)>0$ provided that $\sigma_1\leq\frac1 2$. Let $\sigma_0=2-|\M_b|^2$, when $1<|\M_b^2(r_0)|<2$, one has $1<|\M_b^2(y_1)|<2$ for any $y_1\in[r_0, r_1]$, which can be arrived according to the fact that $\frac{d}{dy_1}|\M_b|^2<0$ for any $r\geq r_0$. Then direct calculations by employing the identity \eqref{22} yield
\begin{align*}
&\quad(w_1k_{b22})'(y_1)=w_1'k_{b22}+w_1k_{b22}'=w_1(2k_{b1}k_{b22}+k_{b22}')-2\sigma_1k_{b22}\\
&=\frac{2|\M_b|^2}{y_1^3(1-M_{b1}^2)^3}(2-|\M_b|^2)\cdot w_1(y_1)-2\sigma_1\frac{1-|\M_b|^2}{y_1^2(1-M_{b1}^2)^2}\\
&\geq\frac{1}{y_1^2(1-M_{b1}^2)^2}\left(\sigma_0\frac{2|\M_b|^2}{y_1(1-M_{b1}^2)}\cdot \left(w_1(r_0)-2\sigma_1\int_{r_0}^{y_1}e^{-\int_{r_0}^{y_1}2k_{b1}(\tau_2)d\tau_2}\right)
-2\sigma_1(1-|\M_b|^2)\right)\\
&\geq\sigma_3, \quad \forall y_1\in[r_0, r_1],
\end{align*}
where the constants $\sigma_2$ and $\sigma_3$ depend only on the background flow and satisfy $0<\sigma_1\leq\sigma_2$.

We now determine $w_2(y_1)$ by setting $$w_2(y_1)=\left(\xi'(r_0)-\frac{w_0}{r_0}\right)w_1(r_0)e^{\int_{r_0}^{y_1}k_{b1}(\tau)d\tau},$$ where $w_0$ is a constant to be selected such that
\be
k_{b22}(r_0)+\left(\xi'(r_0)-\frac{w_0}{r_0}\right)^2>0,
\ee
solving it yields $
w_0<\frac{M_{b1}(r_0)M_{b2}(r_0)-\sqrt{|\M_{b}(r_0)|^2-1}}{1-M_{b1}^2(r_0)} \,\,\mbox{or}\,\, w_0>\frac{M_{b1}(r_0)M_{b2}(r_0)+\sqrt{|\M_{b}(r_0)|^2-1}}{1-M_{b1}^2(r_0)},
$ which agrees with \eqref{112} in Theorem \ref{irro}. Differentiating $w_2$ with respect to $y_1$ yields that
\be\label{l_2}
w_2(y_1)k_{b1}(y_1)-w_2'(y_1)=0,
\ee
and it is easy to compute that
\be\label{5} \begin{split}
&\quad \left(\sqrt{w_1}\p_{y_1}\phi+\frac{w_2}{\sqrt{w_1}}\p_{y_2}\phi\right)(r_0, y_2)\\
&=\frac{\sqrt{w_1(r_0)}}{r_0}(r_0\p_{y_1}\phi+r_0\xi'(r_0)-w_0)\p_{y_2}\phi)(r_0, y_2)\\
&=\frac{\sqrt{w_1(r_0)}}{r_0}g_2(y_2).
\end{split} \ee
If $w_2(y_1)$ fixed, and $\sigma_*$ is chosen as $\sigma_*=\frac1 2 \min\{\sigma_2, \sigma_3\},$ we derive that
\be\label{0} \begin{split}
&\quad\frac1 2\p_{y_1}(w_1k_{22})-\frac1 2w_2\p_{y_2}k_{22}-\p_{y_1}(w_2k_{12})+w_2k_2\\
&=\frac1 2(w_1k_{b22})'(y_1)+\frac1 2\p_{y_1}(w_1(k_{22}-k_{b22}))-\frac1 2w_2\p_{y_2}k_{22}-\p_{y_1}(w_2k_{12})+w_2k_2\\
&\geq\sigma_3-||\p_{y_1}(w_1(k_{22}-k_{b22}))||_{L^{\infty}}-||w_2\p_{y_2}k_{22}||_{L^{\infty}}-||\p_{y_1}(w_2k_{12})||_{L^{\infty}}-||w_2k_2||_{L^{\infty}}\\
&\geq\frac1 2\sigma_3\geq\sigma_*, \quad \forall (y_1, y_2)\in\Omega.
\end{split} \ee
With the aid of the Sobolev embedding $H^3(\Omega)\subset C^{1, \alpha}(\bar{\Omega})$, one can invoke \eqref{estimate} to get
\be\label{1} \begin{split}
w_1k_1-\frac1 2w_1'-w_1\p_{y_2}k_{12}&=w_1k_{b1}-\frac1 2w_1'+w_1(k_1-k_{b1})-w_1\p_{y_2}k_{12}\\
&\geq\sigma_2-||w_1(k_1-k_{b1})||_{L^{\infty}}-||w_1\p_{y_2}k_{12}||_{L^{\infty}}\\
&\geq\frac1 2\sigma_2>0,
\end{split} \ee
here we set $\sigma_1=\sigma_2$. Moreover, owing to $k_{b22}>0$, we deduce that
\be\label{2} \begin{split}
&\quad\frac{w_2^2}{w_1}(r_1, y_2)+k_{22}(r_1, y_2)w_1(r_1)-2k_{12}(r_1, y_2)w_2(r_1)\\
&=\frac{1}{w_1(r_1)}(l_2^2+k_{b22}w_1^2)(r_1)+(k_{22}-k_{b22})(r_1, y_2)w_1(r_1)-2k_{12}(r_1, y_2)|w_2(r_1)|\\
&\geq\frac{1}{w_1(r_1)}(w_2^2+k_{b22}w_1^2)(r_1)-||k_{22}-k_{b22}||_{L^{\infty}}w_1(r_1)-2||k_{12}||_{L^{\infty}}|w_2(r_1)|\\
&\geq\frac{1}{2w_1(r_1)}(w_2^2+k_{b22}w_1^2)(r_1)>0,
\end{split} \ee
and
\be\label{3} \begin{split}
&\quad\frac{w_2^2}{w_1}(r_0, y_2)+k_{22}(r_0, y_2)w_1(r_1)-2k_{12}(r_0, y_2)w_2(r_1)\\
&=w_1(r_0)\left(k_{b22}(r_0)+(\xi'(r_0)-\frac{w_0}{r_0})^2\right)+(k_{22}-k_{b22})(r_0, y_2)w_1(r_0)-2k_{12}(r_0, y_2)|w_2(r_0)|\\
&\geq w_1(r_0)\left(k_{b22}(r_0)+(\xi'(r_0)-\frac{w_0}{r_0})^2\right)-||k_{22}-k_{b22}||_{L^{\infty}}w_1(r_0)-2||k_{12}||_{L^{\infty}}|w_2(r_0)|\\
&\geq \frac{w_1(r_0)}{2}\left(k_{b22}(r_0)+(\xi'(r_0)-\frac{w_0}{r_0})^2\right)>0.
\end{split} \ee
Invoking \eqref{l_2}, one also discovers that
\be\label{4} \begin{split}
&\quad ||k_1w_2-w_2'+w_1k_2-w_1\p_{y_2}k_{22}||_{L^{\infty}}\\
&=||w_2(k_1-k_{b1})-w_1\p_{y_2}k_{22}+w_1k_2||_{L^{\infty}}\\
&\leq C_*\delta_0.
\end{split} \ee
To achieve \eqref{goal}, we integrate \eqref{F} to arrive at
\begin{align*}
&\quad\iint_{\Omega}\left(w_1k_1-\frac1 2w_1'-w_1\p_{y_2}k_{12}\right)(\p_{y_1}\phi)^2dy_1dy_2\\
&\quad+\iint_{\Omega}\left(\frac1 2\p_{y_1}(w_1k_{22})-\frac1 2w_2\p_{y_2}k_{22}-\p_{y_1}(w_2k_{12})+w_2k_2\right)(\p_{y_2}\phi)^2dy_1dy_2\\
&\quad+\iint_{\Omega}(w_2k_1-w_2'+w_1k_2-w_1\p_{y_2}k_{22})\p_{y_1}\phi\p_{y_2}\phi dy_1dy_2\\
&\quad+\int_0^{2\pi}\left(w_1k_{22}-2w_2k_{12}+\frac{w_2^2}{w_1}\right)(\p_{y_2}\phi)^2(r_0, y_2)dy_2\\
&\quad-\frac1 2\int_0^{2\pi}\Bigg(\left(\sqrt{w_1}\p_{y_1}\phi+\frac{w_2}{\sqrt{w_1}}\p_{y_2}\phi\right)^2(r_1, y_2)dy_2\\
&=\iint_{\Omega}\hat{f}(w_1(y_1)\p_{y_1}\phi+w_2(y_1)\p_{y_2}\phi)dy_1dy_2-\frac1 2\int_0^{2\pi}\left(\sqrt{w_1}\p_{y_1}\phi+\frac{w_2}{\sqrt{w_1}}\p_{y_2}\phi\right)^2(r_0, y_2)dy_2\\
&\quad+\int_0^{2\pi}\left(w_1k_{22}-2w_2k_{12}+\frac{w_2^2}{w_1}
\right)(\p_{y_2}\phi)^2(r_1, y_2)dy_2.
\end{align*}
With the aid of \eqref{5}-\eqref{4}, we then obtain that
\begin{equation*} \begin{split}
&\iint_{\Omega}(|\p_{y_1}\phi|^2+|\p_{y_2}\phi|^2)dy_1dy_2+\int_0^{2\pi}\left(\sqrt{w_1}\p_{y_1}\phi+\frac{w_2}{\sqrt{w_1}}\p_{y_2}\phi\right)^2(r_1, y_2)+(\p_{y_2}\phi)^2(r_0, y_2)dy_2\\
&\leq\frac{C_*}{\sigma_*}\left(\iint_{\Omega}|\hat{f}(y_1, y_2)|^2dy_1dy_2+\sum\limits_{j=2}^{3}\int_0^{2\pi}|g_j(y_2)|^2dy_2\right),
\end{split} \end{equation*}
which further implies that
$$
||\phi||_1^2\leq\frac{C_*}{\sigma_*}\left(\iint_{\Omega}\hat{f}^2(y_1, y_2)dy_1dy_2+\sum\limits_{j=2}^{3}\int_0^{2\pi}|g_j(y_2)|^2dy_2\right),
$$
by noting that $\phi(y_1, y_2)=\int_{r_1}^{y_1}\p_{y_1}\phi(\tau, y_2)d\tau+\int_{f(r_1)}^{y_2}g_3(\tau)d\tau$. Then the proof of Lemma \ref{low} have been finished.
\epf
The following lemma gives the higher order estimate.
\bl\label{high}
Under the assumptions of Lemma \ref{low}, the following high order derivatives estimate holds:
\be
||\phi||_4\leq\frac{C_*}{\sigma_*}(||\hat{f}||_3+||g_0||_3+||g_1||_3).
\ee
\el
The proof will largely follow the Lemma 3.6 in \cite{Weng2}, we leave it out here.

\subsection{Proof of Theorem \ref{irro}}

The techniques alluded to in \cite{Weng2} about the proof of Theorem \ref{irro} still work for the isothermal flows in our paper. By means of Galerkin' method, we construct approximate solutions to the nonlinear problem. And the Fourier series is applied to construct the approximate solution to the linearized problem. One may look up the details in \cite{Weng2}(section 3.1.3).\\

\begin{remark}

In fact, the flow may be rotational in general, so it is worthy for us to further investigate the structural stability of smooth transonic isothermal flows with nonzero vorticity. By analogy with the approach taken in \cite{Weng2}, we still employ the deformation-curl decomposition developed in \cite{Weng3, Weng4} so as to handle the elliptic-hyperbolic coupled structure in the steady Euler system. Upon considering that the proof of the structural stability of smooth transonic isothermal flows with nonzero vorticity is similar as in \cite{Weng2} to a large extent, we omit it here, one can refer to section 3.2 of \cite{Weng2}.
\end{remark}

\section{Transonic shock flows in an annulus}\label{transonic}\noindent

The final section devotes to present a precise description of transonic flows with shocks in an annulus $\Omega=\{(x_1, x_2):r_0<r=\sqrt{x_1^2+x_2^2}<r_1\}$.

\begin{definition}(A shock solution of Euler system).\,\,
 We call a piecewise smooth solution $(U_{1}^{\pm},U_{2}^{\pm}, \rho^{\pm})\in C^1(\Omega^{\pm})$ with a jump on the curve $r=\zeta(\theta)$ a shock solution to the system \eqref{radial-euler}, if $(U_{1}^{\pm},U_{2}^{\pm}, \rho^{\pm})$ satisfies the system \eqref{ProblemImm} in $\Omega^{\pm}$ respectively, and the following Rankine-Hugoniot jump conditions:
\begin{align}\label{jump00}
\begin{cases}
[\rho U_1]-\frac{\zeta'(\theta)}{\zeta(\theta)}[\rho U_2]=0,\\
[\rho U_1^2+p]-\frac{\zeta'(\theta)}{\zeta(\theta)}[\rho U_1U_2]=0,\\
[\rho U_1U_2]-\frac{\zeta'(\theta)}{\zeta(\theta)}[\rho U_2^2+p]=0,
\end{cases}
\end{align}
and the entropy condition $[p]>0$, where $[f]$ is defined by $[f(r)]:=f(r)|_{\Omega^+}-f(r)|_{\Omega^-}$.
\end{definition}

There has great ongoing interest in studying the radially symmetric transonic shock problem, with the aim of obtaining more results about its structural and dynamical stability in the mathematical studies of gas dynamics. In \cite{Courant1948, Xin2008a}, the researchers proved the existence and uniqueness of radially symmetric supersonic-subsonic shock solutions without angular velocity in a divergent nozzle. Moreover, in \cite{Xin2008a}, the authors showed that the symmetric transonic shocks were dynamically stable in divergent nozzles and were dynamically unstable in convergent nozzles. The transonic shock problem with given exit pressure was shown to be ill-posed for the potential flow model in \cite{xy05, Xin2008a}. Later on, the symmetric transonic shock solutions of the two-dimensional steady compressible Euler system had been proved structurally stable under the perturbation of the exit pressure and the nozzle wall in \cite{lxy09b}. And the authors further investigated the existence and monotonicity of the axi-symmetric transonic shock by perturbing axi-symmetrically exit pressure in \cite{lxy10a}. Furthermore, the structural stability under axi-symmetric perturbations of the nozzle wall was examined by introducing a modified invertible Lagrangian transformation in \cite{Weng2019}. In the recent researches, the role played by the angular velocity in the structure of steady inviscid compressible flows for Euler system has been elucidated in \cite{Weng1}, and the authors gave a new proof of the existence of smooth transonic spiral flows by studying one side boundary value problem to the steady Euler system in an annulus, which suits their further purpose for the investigation of structural stability of this special transonic flows in \cite{Weng2}. Moreover, \cite{WZ} established the existence of subsonic spiral flows outside a porous body by a variational method.

In our paper, we will prove the existence of radially symmetric transonic solutions with shocks and nonzero angular velocity to \eqref{radial-euler} satisfying suitable boundary conditions both at the inner and outer circle. More precisely, prescribing the pressure at the exit of an annulus $\Omega$, one seeks for smooth functions $(U_{1}^{\pm},U_{2}^{\pm}, \rho^{\pm}, p^{\pm})$ on $\Omega^{\pm}$ respectively with a shock st $r=r_b$ satisfying the physical entropy condition $[p]>0$ and the following Rankine-Hugoniot conditions
\be\label{rankine-hugoniotI}
\begin{cases}
[\rho U_1](r_b)=[\rho U_1^2+p](r_b)=0,\\
[\rho U_1U_2](r_b)=[B](r_b)=0,
\end{cases}
\ee
and the appropriate boundary conditions.

In order to analyze the transonic shock solutions and the wave pattern clearly, it is wise to discuss through dividing the problem into two cases as we did in the previous chapter.

First of all, we consider the transonic shock flows moving from the inner circle to the outer one. In other words, we aim to find smooth functions $(U_{1}^{\pm},U_{2}^{\pm}, \rho^{\pm}, p^{\pm})$ on $\Omega^{\pm}$ respectively with $\Omega^+=\{(x_1,x_2): r_b<r<r_1\}$ and $\Omega^-=\{(x_1,x_2): r_0<r<r_b\}$, supplemented with the boundary conditions:
\begin{eqnarray}\label{bd3}
\rho(r_0)=\rho_{0}>0,\ U_{1}(r_0)=U_{10}>0, \ U_{2}(r_0)=U_{20}(\neq 0), \ p(r_1)= p_{ex}.
\end{eqnarray}

For later use, we readily define the functions
\be
\psi_1(r)=\frac{A-(U_2^-)^2}{A}\cdot\frac{(U_1^-)^2}{A},
\ee
\begin{align}\label{f_2}
\psi_2(\rho):=p_{ex}\cdot\ln\left(\frac{p_{ex}}{A}\right)\rho-(B_0-\dfr12
\dfr{m_2^2}{r^2})\rho^2+\dfr{m_1^2}{2r^2}.
\end{align}

\bt\label{ProbIII}
{\it For a given incoming flow which is supersonic in the $r$-direction and has a nonzero angular velocity at the entrance, there exist two constants $p_1<p_0$ such that, for any $p_{ex}\in(p_1,p_0)$, there exist a unique piecewise smooth solution to \eqref{radial-euler} with the boundary conditions \eqref{bd3} in $\Omega$ with a shock located at $r\equiv r_b$ and satisfying $U_1^-(r)>A$ for $r\in(r_0,r_b)$, $U_1^+(r)<A$ for $r\in(r_b,r_1)$. Furthermore, the solution has the following properties.
\begin{enumerate}[(A)]
  \item The shock position $r=r_b$ increases as the exit pressure $p_{ex}$ decreases. In addition, the shock position $r_b$ tends to $r_1$ if $p_{ex}$ goes to $p_1$ while $r_b$ approaches to $r_0$ if $p_{ex}$ goes to $p_0$.
  \item The flow patterns in $\Omega^+$ can be classified in terms of the boundary conditions as follows:
  \begin{enumerate}[1.]
  \item If $r_0<r_1\leq r_{**}$ and $p_{ex}\in(p_1,p_0)$, then
  \begin{description}
  \item{Subcase 1.1.} for $\psi_2(\rho_s)>0$, there exists a supersonic-supersonic shock and the flow is supersonic in $\Omega^+$.
  \item{Subcase 1.2.} for $\psi_2(\rho_s)<0$, there exists a supersonic-supersonic shock and the flow changes smoothly from supersonic to subsonic in $\Omega^+$.
  \item{Subcase 1.3.} for $\psi_2(\rho_s)=0$, there exists a supersonic-supersonic shock and the flow is supersonic in $\Omega^+$ but degenerates to a sonic state at the exit.
  \end{description}
  \item If $r_0<r_{**}<r_1$, then there exists a $p_{**}\in(p_1,p_0)$ such that,
  \begin{description}
  \item{Subcase 2.1.} for $p_{ex}\in(p_{**},p_0)$ with $\psi_2(\rho_s)>0$, a supersonic-supersonic shock exists and the flow is supersonic in $\Omega^+$.
  \item{Subcase 2.2.} for $p_{ex}\in(p_{**}p_0)$ with $\psi_2(\rho_s)<0$, a supersonic-supersonic shock exists and the flow changes smoothly from supersonic to subsonic in $\Omega^+$.
  \item{Subcase 2.3.} for $p_{ex}\in(p_{**},p_0)$ with $\psi_2(\rho_s)=0$, a supersonic-supersonic shock exists and the flow is supersonic in $\Omega^+$ but degenerates to a sonic state at the exit.
  \item{Subcase 2.4.} for $p_{ex}\in(p_1,p_{**})$, a supersonic-subsonic shock exists and the flow is subsonic in $\Omega^+$.
  \item{Subcase 2.5.} for $p_{ex}=p_{**}$, there exists a supersonic-sonic shock and the flow is subsonic in $\Omega^+$ but degenerates to a sonic state at the shock position.
  \end{description}
  \item If $ r_1>r_0\geq r_{**}$, then for any $p_{ex}\in(p_1,p_0)$, there exists a supersonic-subsonic shock and the flow is uniformly subsonic in $\Omega^+$.
\end{enumerate}
\end{enumerate}
}
\et

\bpf
To start with, we see from \eqref{alg1} that
\begin{align}
\begin{cases}
\rho^- U_1^-=\dfr{m_1}{r},\ m_1= r_0 \rho_{0} U_{10},\\
U_2^-=\dfr{m_2}{r},\ m_2= r_0 U_{20},\\
B^-\equiv B_0, B_0=\dfrac{1}{2}(U_{10}^2+U_{20}^2)+A \ln\rho_0
\end{cases}
\end{align}
in $\Omega^-$. In addition, one has
\begin{align*}\label{Omgp}
(\rho^+ U_1^+)(r)=\dfr{m_1}{r},\ \ U_2^+(r)=\dfr{m_2}{r}
\end{align*}
in $\Omega^+$, with the aid of \eqref{rankine-hugoniotI} which implies that $m_1$, $m_2$ and are unchanged across the shock.

{\textbf{Step 1.}}\quad First, we set out to determine the relationship between the shock position $r=r_b$ and the exit pressure $p(r_1)=p_{ex}$.

In view that on the exit $r=r_1$, it holds that
\be
p_{ex}=A\rho^+(r_1),
\ee
so it is clear that
\be\label{dp}
\dfr{d p_{ex}}{dr_b}=A\dfr{d\rho^+(r_1)}{dr_b}.
\ee
Moreover, the Bernoulli's function can be written as
\be\no
B^+(r_1)=\dfr12(\dfr{m_1^2}{(\rho^+(r_1))^2 r_1^2}+\dfr{m_2^2}{r_1^2})+A\ln\rho^+(r_1),
\ee
and
\be\no
B^+(r_b)=\dfr12(\dfr{m_1^2}{(\rho^+(r_b))^2 r_b^2}+\dfr{m_2^2}{r_b^2})+A\ln\rho^+(r_b)
\ee
on $ r=r_1$ and the shock $r=r_b$ respectively. Then the employment of the Bernoulli's laws gives
\be
\dfr12(\dfr{m_1^2}{(\rho^+(r_1))^2r_1^2}+\dfr{m_2^2}{r_1^2})+A\ln\rho^+(r_1)=\dfr12(\dfr{m_1^2}{(\rho^+(r_b))^2 r_b^2}+\dfr{m_2^2}{r_b^2})+A\ln\rho^+(r_b).
\ee
Differentiating this with respect to $r_b$ yields
\be
\frac{1}{\rho^+(r_1)}\left(\frac{m_1^2}{(\rho^+(r_1))^2r_1^2}-A\right)\frac{d\rho^+(r_1)}{dr_b}
=\frac{1}{\rho^+(r_b)}\left(\frac{m_1^2}{(\rho^+(r_b))^2r_b^2}-A\right)\frac{d\rho^+(r_b)}{dr_b}+\frac{m_1^2}{(\rho^+(r_b))^2r_b^3}+\frac{m_2^2}{r_b^3}.
\ee
In addition, the R-H conditions \eqref{rankine-hugoniotI} implies that
\be
\frac{m_1^2}{r_b^2\rho^+(r_b)}+A\rho^+(r_b)=\frac{m_1^2}{r_b^2\rho^-(r_b)}+A\rho^-(r_b).
\ee
Continue to differentiate it with respect to $r_b$ indicts that
$$
\frac{2m_1^2}{r_b^3\rho^+(r_b)}+\left(\frac{m_1^2}{r_b^2(\rho^+(r_b))^2}-A\right)\frac{d\rho^+(r_b)}{dr_b}
=\frac{2m_1^2}{r_b^3\rho^-(r_b)}+\left(\frac{m_1^2}{r_b^2(\rho^-(r_b))^2}-A\right)\frac{d\rho^-(r_b)}{dr_b},
$$
combining this with \eqref{daoshu}, one has
\be
\frac{1}{\rho^+(r_b)}\left(\frac{m_1^2}{(\rho^+(r_b))^2r_b^2}-A\right)\frac{d\rho^+(r_b)}{dr_b}
=-\frac{2m_1^2}{r_b^3}\frac{1}{(\rho^+(r_b))^2}-\frac{m_2^2}{r_b^3}\frac{\rho^-(r_b)}{\rho^+(r_b)}+\frac{m_1^2}{r_b^3\rho^-(r_b)}\frac{1}{\rho^+(r_b)}.
\ee
Herein, we derive that
\begin{align}
&\quad\frac{1}{\rho^+(r_1)}\left(\frac{m_1^2}{(\rho^+(r_1))^2r_1^2}-A\right)\frac{d\rho^+(r_1)}{dr_b}\\
&=-\frac{2m_1^2}{r_b^3}\frac{1}{(\rho^+(r_b))^2}-\frac{m_2^2}{r_b^3}\frac{\rho^-(r_b)}{\rho^+(r_b)}+\frac{m_1^2}{r_b^3\rho^-(r_b)}\frac{1}{\rho^+(r_b)}
+\frac{m_1^2}{(\rho^+(r_b))^2r_b^3}+\frac{m_2^2}{r_b^3}\no\\
&=\frac{m_1^2}{r_b^3(\rho^+(r_b))^2}\left(\frac{\rho^+(r_b)}{\rho^-(r_b)}-1\right)+\frac{m_2^2}{r_b^3}\left(1-\frac{\rho^-(r_b)}{\rho^+(r_b)}\right)>0,
\end{align}
with the help of $\rho^-(r_b)<\rho^+(r_b)$. Thus it holds that
\be\label{d rho1}
\frac{d\rho^+(r_1)}{dr_b}<0,
\ee
owing to $\frac{m_1^2}{(\rho^+(r_1))^2r_1^2}-A<0$. Then one can conclude from \eqref{dp} and \eqref{d rho1} that
\be
\dfr{d p_{ex}}{dr_b}<0.
\ee

We have so far shown the following result.
\bp\label{prop36}
For any given incoming flow $U_1^-(r_0), U_2^-(r_0), \rho^-(r_0), p^-(r_0)$ satisfying $U_1^-(r_0)>A$, there exist two constant $p_1<p_0$ such that, for any $p_{ex}\in(p_1, p_0)$, problem \eqref{radial-euler} and \eqref{bd3} has a unique solution with a shock located at $r=r_b$, satisfying $U_1^-(r)>A$ for $r\in(r_0, r_b)$, $U_1^-(r)<A$ for $r\in(r_b, r_1)$. In addition, the shock position $r_b$ increases as the exit pressure $p_{ex}$ decreases. Furthermore, the shock position $r_b$ approaches to $r_1$ if $p_{ex}$ tends to $p_1$ while $r_b$ approaches to $r_0$ if $p_{ex}$ tends to $p_0$.
\ep
{\textbf{Step 2.}}\quad It remains to clarify the flow state just behind the shock. In other words, the next task is to determine the sign of
\be\no
(U_{1}^+)^2+(U_{2}^+)^2 -A
\ee
at $r=r_b$. It follows from the first two equations in \eqref{rankine-hugoniotI} that $[U_1+\frac{A}{U_1}]=0$, which means that $U_1^+(r_b)=\frac{A}{U_1^-(r_b)}$, thus one has
\begin{align}\nonumber
\quad(U_{1}^+)^2+(U_{2}^+)^2-A&=\frac{A^2}{(U_1^-)^2}+(U_{2}^-)^2-A\\\no
&=\frac{A}{(U_1^-)^2}\left(A+\frac{(U_1^-)^2(U_2^-)^2}{A}-(U_1^-)^2\right)\\\no
&=\frac{A}{(U_1^-)^2}\Bigg(\left((U_2^-)^2-A\right)\frac{(U_1^-)^2}{A}+A\Bigg).
\end{align}
Let us for simplicity to define
\begin{align*}
x(r)=\frac{(U_1^-)^2}{A}=(M_1^-(r))^2,  \quad\mbox{and}   \quad h(x)=\left((U_2^-)^2-A\right)\frac{(U_1^-)^2}{A}+A,
\end{align*}
then one may write
\begin{align*}
\psi_1(r)=\frac{(A-(U_2^-)^2)}{A}x.
\end{align*}
By virtue of
\be\label{ainfty}
\frac{d}{dr}((U_2^-)(r))^2<0,\ \ \lim_{r\rightarrow +\infty}((U_2^-)(r))^2=0,
\ee
there exists a unique $r_*=\sqrt{\frac{m_2^2}{A}}$ such that $((U_2^-)(r_*))^2=A$. Since $((U_2^-)(r))^2<A$ for any $r> r_*$ and decreases as $r$ increases, $x(r)>1$ and $x'(r)>0$ for any $r>r_0$, one has that $\psi_1'(r)>0$ for all $r>r_*$. We also note that $(M_1^2)'\geq\frac1 rM_1^2$ by \eqref{daoshu}, which implies that $M_1^2(r)\geq\frac{r}{r_0}M_1^2(r_0)$. This, together with \eqref{ainfty} yields
\be\nonumber
\lim_{r\rightarrow +\infty}\psi_1(r)=+\infty.
\ee
Since $\psi_1(r_*)=0$ and $\displaystyle\lim_{r\to \infty} \psi_1(r)=\infty$, there exists a unique $r_{**}>r_*$ such that $\psi_1(r_{**})=1$.

For the sake of determining the state behind the shock, let us restrict attention to the judgment of the sign of $h(x(r_b))$, by noticing that $(U_1^-)^2+(U_2^-)^2-A=h(x)\cdot x$.

To proceed it, we compute the root of $h(x)=0$:
$$
x_0=\frac{A}{A-(U_2^-)^2}.
$$
Henceforth, the sign of $h(x(r_b))$ can be analyzed in detail as follows.
\begin{enumerate}[(1)]
  \item $(U_2^-)^2=A$. Then $h(x(r_b))\equiv A>0$ for all $x(r_b)\geq0$, which indicates that the flow behind the shock is supersonic.
  \item $(U_2^-)^2>A$. Then $h(x(r_b))>h(0)=A>0$ for all $x(r_b)\geq0$ by noting that $h(x(r))$ increases monotonically with respect to $r$, thus the flow behind the shock is supersonic as well.
  \item $(U_2^-)^2<A$. Then $h(x(r_b))>0$ for $0\leq x(r_b)<\frac{A}{A-(U_2^-)^2}$, while $h(x(r_b))<0$ for $x(r_b)>\frac{A}{A-(U_2^-)^2}$. Thus, it holds that
      \begin{enumerate}[(i)]
      \item if $\psi_1(r_b)<1$, then the state behind the shock is supersonic;
      \item if $\psi_1(r_b)>1$, then the state behind the shock is subsonic;
      \item if $\psi_1(r_b)=1$, then the state behind the shock is sonic.
      \end{enumerate}
\end{enumerate}

Take notice that the position $r_b$ of the shock is unknown, the sign of $h(x(r_b))$ can not even be determined. To solve this, it is wise to use the properties of $(U_2^-)^2$, $\psi_1$ and the monotonicity between the shock position and the exit pressure in three different cases:

{\bf Case a}: $r_0<r_*$ (i.e. $(U_{20}^-)^2>A$).  The flow is supersonic (subsonic) behind the shock if the shock occurs at $r_0<r_b<r_{**}$ ($r_b>r_{**}$), and is sonic behind the shock if $r_b=r_{**}$. Since the exit pressure $p_{ex}$ depends monotonically on the shock position $r=r_b$, there exists an interval $(p_1,p_0)$ such that if the exit pressure $p_{ex}\in(p_1,p_0)$, there exists a unique shock solution with a shock at $r_b\in(r_0,r_1)$.

If $r_1\leq r_{**}$, one obtains a supersonic-supersonic shock. Note that the flow at downstream can change smoothly from supersonic to subsonic after crossing the shock.

If $r_1> r_{**}$, then there exists a $p_{**}\in(p_1,p_0)$ corresponding to $r_{**}$ such that
\begin{enumerate}[(i)]
  \item if $p_{ex}\in(p_1, p_{**})$, one has a supersonic-subsonic shock located at $r_b\in(r_{**}, r_1)$;
  \item if $p_{ex}\in(p_{**}, p_0)$, one gets a supersonic-supersonic shock located at $r_b\in(r_0, r_{**})$ and the flow at downstream can change smoothly from supersonic to subsonic after crossing the shock.
  \item if $p_{ex}=p_{**}$, one obtains a supersonic-sonic shock with a shock and the sonic circle both located at $r_b=r_{**}$.
\end{enumerate}

{\bf Case b}:  $r_*\leq r_0<r_{**}$.  The flow is supersonic behind the shock if the shock occurs at $r_0<r_b<r_{**}$, and is subsonic (sonic) behind the shock if the shock occurs at $r_b>r_{**} (r_b=r_{**})$. Therefore there exists an interval $(p_1, p_0)$ such that if the exit pressure $p_{ex}\in(p_1,p_0)$, there exists a unique shock solution with a shock at $r_b\in(r_0,r_1)$.

  If $r_1\leq r_{**}$, then the shock is supersonic-supersonic, the flow at downstream can change smoothly from supersonic to subsonic crossing the shock.

  If $r_1> r_{**}$, then there exists a $p_{**}\in(p_1,p_0)$ corresponding to $r_{**}$ such that
  \begin{enumerate}[(i)]
  \item if $p_{ex}\in(p_1, p_{**})$, one gets a supersonic-subsonic shock with the shock located at $r_b\in(r_{**},r_1)$;
  \item if $p_{ex}\in(p_{**}, p_0)$, one has a supersonic-supersonic shock with the shock located at $r_b\in(r_0,r_{**})$, the flow at downstream can change smoothly from supersonic to subsonic crossing the shock;
  \item if $p_{ex}=p_{**}$, one obtains a supersonic-sonic shock with the shock and the sonic line both located at $r_b=r_{**}$.
  \end{enumerate}

{\bf Case c}: $r_0\geq r_{**}>r_*$. Thanks to $\psi_1(r_0)\geq \psi_1(r_{**})=1$, $\psi_1(r)>1$ for any $r>r_0$. The flow is subsonic behind the shock if the shock occurs at $r_b>r_0$. Therefore, there exists an interval $(p_1,p_0)$ such that, for any $p_{ex}\in(p_1,p_0)$, there exists a unique supersonic-subsonic shock solution with the shock located at $r_b\in(r_0,r_1)$.

Up to now, we have proved that, for the solutions declared in Theorem \ref{ProbIII}, the flow state just behind the shock can be divided into three cases summarized as follows.
\bp\label{prop37}
For the solutions given in Proposition \ref{prop36}, the state of flow just behind the shock can be classified as:
\begin{enumerate}[1.]
  \item If $r_0<r_1\leq r_{**}$, then $|{\bf M}^+|^2(r_b)>1$ for any $p_{ex}\in(p_1, p_0)$.
  \item If $r_0<r_{**}<r_1$, there exists a $p_{**}\in(p_1,p_0)$ such that $|{\bf M}^+|^2(r_b)>1$ for any $p_{ex}\in(p_{**}, p_0)$, $|{\bf M}^+|^2(r_b)<1$ for any $p_{ex}\in(p_1,p_{**})$ and $|{\bf M}^+|^2(r_b)=1$ for $p_{ex}=p_{**}$.
  \item If $r_0\geq r_{**}$, then $|{\bf M}^+|^2(r_b)<1$ for any $p_{ex}\in(p_1,p_0)$.
  \end{enumerate}
\ep

{\textbf{Step 3.}}\quad Finally, it is ought for us to demonstrate the flow state at the exit by determining the sign of
\be\no
(U_{1}^+)^2+(U_{2}^+)^2 -A
\ee
at $r=r_1$.

Upon concerning \eqref{f_2}, it is easy to get that
$$
\psi_2(\rho)=p_{ex}(\ln p_{ex}-\ln A)-2\left(B_0-\frac1 2\frac{m_2^2}{r^2}\right)\rho,
$$
and
$$
\psi_2(\rho^+(r_1))=0,
$$
which tells that $\psi_2(\rho)$ increases on $(0, \rho_{\sharp})$ while decreases on $(\rho_{\sharp}, \infty)$ as $\rho$ increases, where
$$
\rho_{\sharp}=\frac{p_{ex}(\ln p_{ex}-\ln A)}{2(B_0-\frac1 2\frac{m_2^2}{r_1^2})}.
$$
In light of $\psi_2(0)>0$, one can say that $\rho<\rho^+(r_1)$ if and only if $\psi_2(\rho)>0$, while $\rho>\rho^+(r_1)$ if and only if $\psi_2(\rho)<0$.

By noticing that that the flow is supersonic at $r=r_1$ if and only if $B_0-A\ln{\rho^+(r_1)}>\frac A 2$, that is to say,
\be\label{rhosharp}
\rho^+(r_1)<\exp\left(\frac{B_0}{A}-\frac1 2\right)=\rho_s \,(\mbox{defined in}\, \eqref{rhoc}),
\ee
in agreement with $\psi_2(\rho_s)>0$. In turn, the flow is subsonic at $r=r_1$ if and only if $\rho^+(r_1)<\rho_s$, which coincides with $\psi_2(\rho_s)<0$.

To sum up, we have obtained the following results.
\bp\label{prop38}
For the solutions given in Proposition \ref{prop36}, the state of flow at the exit can be classified as
\begin{enumerate}
  \item If $\psi_2(\rho_s)>0$, then $|\M^+|^2(r_1)>1$.
  \item If $\psi_2(\rho_s)<0$, then $|\M^+|^2(r_1)<1$.
  \item If $\psi_2(\rho_s)=0$, then $|\M^+|^2(r_1)=1$.
\end{enumerate}
\ep

Combining Proposition \ref{prop36}, \ref{prop37} and \ref{prop38} obtained above, we completed the proof of Theorem \ref{ProbIII}. \epf

\begin{remark}
The analysis for the transonic shock wave patterns when a supersonic flow moves from the outer to the inner circle is similar as the the case where the flow moves from the inner to the outer circle. Slightly different from the latter one, there exists a limiting circle $r=r^{\sharp}$, to which the acceleration of the solution established in Theorem \ref{subsec2} will blow up as $r$ tends $r^{\sharp}$. For this reason we shall henceforth introduce $\rho_e=\exp\left(\frac{K(r_0)}{A}-\frac1 2\right)$ to compare $r_0$ and $r^{\sharp}$. According to the different situations of $\psi_2(\rho_e)$ and the exit pressure $p_{ex}$, the different states of solution can be revealed, and we also give the corresponding flow patterns in terms of the boundary conditions similarly as in Theorem \ref{ProbIII}.
\end{remark}

\end{document}